\documentclass{amsart}
\usepackage{amssymb}
\usepackage{graphicx}
\usepackage{epstopdf}
\input {epsf.sty}
\newtheorem{thm}{Theorem}[section]

\newtheorem{lem}[thm]{Lemma}
\newtheorem{prop}[thm]{Proposition}
\theoremstyle{definition}
\newtheorem{defn}[thm]{Definition}

\numberwithin{equation}{section}

\newcommand{\R}{\mathbb R}

\newcommand{\To}{\longrightarrow}

\newcommand{\C}{\mathbb{C}}

\newcommand{\inv}{^{-1}}

\newcommand{\x}{\times}

\newcommand{\K}{\overline{\mathbb K}}

\begin{document}
\title[Totally nonnegative cells in $G/P$]
{Closure relations for totally nonnegative cells in $G/P$.}
\author{K. Rietsch}%
\address{Department of Mathematics,
            King's College London,
            Strand, London
            WC2R 2LS}%
\address{Department of Pure Mathematics, University of Waterloo,
Canada}
\email{rietsch@mth.kcl.ac.uk}%
\thanks{
The author is supported by a Royal Society Dorothy
Hodgkin Research Fellowship.}%
%
\begin{abstract}
The totally nonnegative part of a partial flag variety $G/P$ 
has been shown in \cite{Rie:CelDec,Rie:PhD} to be a union
of semi-algbraic cells.  
We show that the closure of a cell is a union of cells and 
give a combinatorial description of the closure relations. 
The totally nonnegative cells are defined by intersecting
$(G/P)_{\ge 0}$ with a certain stratification of $G/P$
defined by Lusztig \cite{Lus:PartFlag}.
We also verify the same closure relations for these strata. 
\end{abstract}
\maketitle

\section{Introduction} 
For a reductive algebraic group over $\C$ split over $\R$ 
with fixed choice of Chevalley generators in the Lie algebra, there
is a well defined notion of positive, or $\R_{> 0}$-valued, points
due to Lusztig \cite{Lus:TotPos94}. 
In the case of $GL_n$ with the standard choices the resulting 
``$GL_n(\R_{>0})$'' recovers the classical notion
of totally positive matrices,
that is matrices all of whose minors are in $\R_{>0}$.  
For general $G$ the set $G(\R_{>0})$, or $G_{>0}$
as we will denote it, is therefore called  
the totally positive part of $G$.
The closure $G_{\ge 0}$ of $G_{>0}$ (in the real topology)
is called the 
totally nonnegative part of $G$.

These notions extend in a natural way to 
flag varieties $G/P$, \cite{Lus:TotPos94,Lus:IntroTotPos}. 
That is, one has a notion of 
$(G/P)_{>0}$ -- this is a semi-algebraic
subset of the real points in $G/P$ --
and of $(G/P)_{\ge 0}$, the closure 
of $(G/P)_{>0}$. 
Now recall that $G/B$ has a decomposition into smooth strata 
$\mathcal R_{v,w}$, obtained as
intersections of opposed Bruhat cells
and indexed by pairs $v,w$ in the Weyl group
with $v\le w$.  
In \cite{Lus:PartFlag} Lusztig defined an
analogous decomposition of $G/P$ into
smooth strata $\mathcal P_{x,u,w}$.
These decompositions intersected with the
$(G/P)_{\ge 0}$ give cell decompositions 
of the totally nonnegative parts of the $G/P$,
\cite{Rie:CelDec,Rie:PhD}. 
We call the components of this decomposition
of $(G/P)_{\ge 0}$ 
the totally nonnegative cells
in $G/P$. Note that there is one open totally nonnegative
cell in $G/P$, namely $(G/P)_{>0}$
itself. 

It was proved in \cite{Lus:IntroTotPos} that 
$(G/B)_{\ge 0}$ is contractible, and the same holds for $(G/P)_{\ge 0}$
by the same proof. Also Fomin and Shapiro \cite{FoSha:StratSpaces} 
studied links of totally nonnegative cells inside a big cell of $SL_n/B$, 
in particular showing them to be contractible. Beyond these special cases, however, not 
much is known
about the closures of the individual cells
or how the cells are glued together.   

In this paper we prove that the closure of a totally 
nonnegative cell is a union of totally nonnegative cells
and describe the closure relations in terms of the 
Weyl group. 
In the full flag variety case we show that  
$\mathcal R^{>0}_{v',w'}\subseteq
\overline{\mathcal R^{>0}_{v,w}}$, 
whenever $v\le v'\le w'\le w$, see Theorem~\ref{t:B}.   
This theorem is then generalized to $G/P$ (see Theorem~\ref{t:P}), with
the combinatorial description of the cells and 
their closure relations given in  Section~\ref{s:P}.

\begin{figure}
\begin{center}
\leavevmode
 \[ \includegraphics{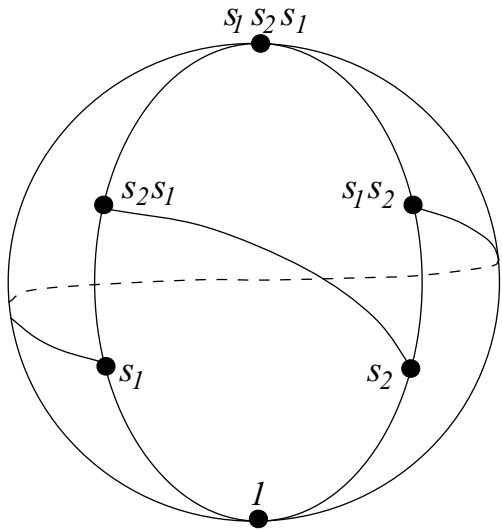}\]
\caption{$(SL_3/B)_{\ge 0}$ with its cell decomposition}
\end{center}
\label{f:SL3}
\end{figure}

A main difficulty in proving such results lies in how to find  
a totally nonnegative cell $\mathcal R_{v,w}^{>0}$ 
inside the closure of another. While much 
detailed information 
is available about the individual cells, such as parameterizations, 
explicit defining equalities/inequalities (see \cite{MarRie:ansatz}), none of 
these results readily extend to the closures of the cells.  The central idea
which allows us to relate this problem to 
an easier special case is contained in Lemma~4.3.

The combinatorial properties of our poset describing the closure relations between
totally nonnegative cells  were recently investigated by Williams \cite{Williams:poset}. 
Her results suggest, and she conjectures, that  
$(G/P)_{\ge 0}$ is a regular CW complex homeomorphic
to a closed ball.  
For an illustration of the totally nonnegative part with its cell 
decomposition in the case of $SL_3/B$
see Figure~\ref{f:SL3}.

In the last section we
verify the same closure  
relations as among the totally nonnegative cells, for the strata 
$\mathcal P_{x,u,w}$ of $G/P$. In this case $G/P$ is either taken 
again over $\R$, 
or over an algebraically closed
field $\mathbb K$ (and Zariski topology). 
The closure relations among the totally nonnegative
cells could in retrospect be viewed as an $\R_{>0}$-valued 
analogue of Proposition~\ref{p:P}, although the results
involving positivity are more difficult to prove. 
  
Lusztig's stratification of $G/P$ (over $\C$) has recently been reinterpreted
by Goodearl and Yakimov \cite{GoodYak:PoissonG/P} in a Poisson geometric
setting. Namely the strata arise as 
torus-orbits of symplectic leaves for a certain natural Poisson structure on $G/P$. 
This paper also independently gives closure relations among these strata 
\cite[Theorem~1.8]{GoodYak:PoissonG/P} which look quite different
from our Proposition~\ref{p:P}. The 
combinatorial equivalence of their description of the poset structure with ours
was recently proved by Xuhua He \cite{Williams:PersonalCommunication}
using arguments from \cite{He:GStablePieces}.

\vskip .2cm
\noindent{\it Acknowledgements.} I am grateful to Lauren Williams 
for prompting me to think about the closure 
relations between the totally nonnegative cells and
for useful discussions.

\section{Preliminaries}\label{s:prelim}
\subsection{}
We recall some basic notation and results from algebraic groups,
see e.g. \cite{Springer:AlgGroupBook}. Suppose 
$\mathbb K$ is an algebraically
closed field, $\mathbb K=\K$, or $\mathbb K=\R$.  
Let $G$ be a semisimple linear algebraic group over 
$\K$ split over $\mathbb K$. We identify $G$ and any related spaces with their
$\mathbb K$-valued points. If $\mathbb K=\R$  we consider 
them with their real topology (as real manifolds or subsets thereof),
otherwise we consider their Zariski topology.   

Let $T$ be a split torus and $B^+$ and $B^-$ opposite Borel 
subgroups containing $T$. The unipotent radicals of $B^+$ and 
$B^-$ are denoted $U^+$ and $U^-$, respectively. 
Let $\{\alpha_i\ | \ i\in I\}$ be the set of simple roots
associated to $B^+$ and $\{\alpha_i^\vee\ |\ i\in I\}$ the 
corresponding coroots. Then we have the simple root subgroups 
$U^+_{\alpha_i}\subseteq U^+$ and $U^-_{\alpha_i}\subseteq U^-$.
Furthermore assume we are given homomorphisms
\begin{equation*}
\phi_i: SL_2(\mathbb K)\to G, \qquad i\in I,
\end{equation*}  
such that 
\begin{equation*}
\phi_i\left(\begin{pmatrix}t &0\\ 0& t\inv 
\end{pmatrix}\right)=\alpha_i^\vee(t), \qquad t\in \mathbb K^*,
\end{equation*}
and such that 
\begin{equation*}
\phi_i\left(\begin{pmatrix}1 & m\\ 0&1\end{pmatrix}\right):=x_i(m), 
\quad 
\phi_i\left(\begin{pmatrix}1 & 0\\ m&1\end{pmatrix}\right):=y_i(m), 
\end{equation*}
define isomorphisms $x_i:\mathbb K\to U^+_{\alpha_i}$ and $y_i:
\mathbb K\to U^-_{\alpha_i}$.

Let $W=N_G(T)/T$ be the Weyl group of $G$. For $i\in I$ the elements
\begin{equation*}
\dot s_i=x_i(-1)y_i(1)x_i(-1)
\end{equation*} 
represent the simple reflections $s_i\in W$. If $w=s_{i_1}\dotsc s_{i_m}$ 
is a reduced expression for $w$ then we write $\ell(w)=m$ for the length
of $w$. It is also known that the representative
\begin{equation*}
\dot w=\dot s_{i_1}\dotsc \dot s_{i_m}
\end{equation*}   
of $w$ is well defined, independent of the choice of reduced expression. 
Inside $W$ we have a longest element which is denoted by $w_0$.

\subsection{}
Let $J\subseteq I$. The parabolic subgroup 
$W_J\subseteq W$ is the subgroup generated by all of the
$s_j$ with $j\in J$. Let $w_J$ denote the longest element in 
$W_J$.  We also consider the set $W^J$
of minimal coset representatives for $W/W_J$, and the
set $W^J_{max}=W^Jw_J$ of maximal coset representatives.

The parabolic subgroup $W_J$ of $W$
corresponds to a parabolic subgroup $P_J$ 
in $G$ containing $B^+$. Namely, 
$P_J$ is the subgroup of $G$ generated by $B^+$ and the
elements $\dot w$ for $w\in W_J$.
Let  $\mathcal P^J$ be the set of parabolic subgroups $P$  
conjugate to $P_J$. This is a homogeneous space for the conjugation 
action of $G$ and can be identified with the 
partial flag variety $G/P_J$ via
 $$
G/P_J\overset\sim\To \mathcal P^J\ :\ gP_J\mapsto gP_J g\inv.
 $$
In the case $J=\emptyset$ we are identifying the full flag variety $G/B^+$
with the variety $\mathcal B$ of Borel subgroups in $G$. We have the
usual projection from the full flag variety to any 
partial flag variety which takes the form 
$\pi=\pi^J:\mathcal B\to\mathcal P^{J}$, where $\pi(B)$ is the 
unique parabolic subgroup of type $J$ containing $B$. 

The conjugate of a parabolic subgroup $P$ by 
an element $g\in G$ will be denoted by
$g\cdot P:=gPg\inv$.  

\subsection{}
Recall the Bruhat decomposition for the full flag variety,
$$
\mathcal B=\bigsqcup_{w\in W} B^+\dot w\cdot B^+,
$$
and the Bruhat order $\le$ on $W$. 
The Bruhat cell $B^+\dot w\cdot B^+$ is isomorphic
to $\mathbb K^{\ell(w)}$. And the Bruhat order has the property
$$
v\le w\ \  \iff \ \  
B^+\dot v\cdot B^+\subseteq\overline {B^+\dot w\cdot B^+},
$$
for $v,w\in W$.

It is a well known consequence of Bruhat decomposition 
that $\mathcal B\x\mathcal B$ is the union of the 
$G$-orbits $\mathcal O(w)=G\cdot (B^+,\dot w\cdot B^+)$,
with $G$ acting
diagonally. Therefore to any pair $(B_1,B_2)$ of Borel 
subgroups one can associate a unique $w\in W$ such that
 \begin{equation*}
(B_1,B_2)=(g\cdot B^+,g\dot w \cdot B^+)
\end{equation*} 
for some $g\in G$. We write  
$$
B_1\overset w\to B_2
$$ 
in this case and call $w$ the relative position of $B_1$ and $B_2$.

\subsection{}
Finally, let us consider the two opposite Bruhat 
decompositions 
\begin{equation*}
\mathcal B=\bigsqcup_{w\in W} B^+\dot w\cdot B^+=\bigsqcup_{v\in W}
B^-\dot v\cdot B^+.
\end{equation*}
Note that $B^-\dot v\cdot B^+\cong\mathbb K^{\ell(w_0)-\ell(v)}$. The 
closure relations for these opposite Bruhat cells are 
given by $B^-\dot v'\cdot B^+\subset 
\overline{B^-\dot v\cdot B^+}$ if and only if $v\le v'$. 
We define
\begin{equation*}
\mathcal R_{v,w}:=B^+\dot w\cdot B^+\cap B^-\dot v\cdot B^+,
\end{equation*}
the intersection of opposed Bruhat cells. This intersection is empty
unless $v\le w$, in which case it is smooth of dimension 
$\ell(w)-\ell(v)$, see \cite{KaLus:Hecke2,Lus:PartFlag}.

\section{Total Positivity for $G$ and $\mathcal B$}
Let $\mathbb K=\R$. The totally nonnegative part $G_{\ge 0}$ of $G$ is defined
by Lusztig \cite{Lus:TotPos94} to be the semigroup
inside $G$ generated by the sets
\begin{align*}
&\{x_i(t)\ |\ t\in \R_{>0}, i\in I\},\\
&\{y_i(t)\ |\ t\in \R_{>0}, i\in I\}, \text{ and}\\
&T_{>0}:=\{ t\in T\ |\ \text{$\chi(t)>0$ all $\chi\in X^*(T)$}\}.
\end{align*}
When $G=SL_n(\R)$ then by a Theorem of A.~Whitney's 
this definition agrees with the classical notion 
of totally nonnegative matrices inside $SL_n(\R)$, 
that is those matrices all of whose minors are 
nonnegative.

\subsection{}\label{s:totpos}
We recall some basic facts about 
total positivity for $G$ from \cite{Lus:TotPos94}.  
Let 
$U^+_{\ge 0}:=G_{\ge 0}\cap U^+$  and $U^-_{\ge 0}:=G_{\ge 0}\cap U^-$.
For $w\in W$ and $s_{i_1}\dotsc s_{i_m}=w$ a reduced expression define
\begin{eqnarray*}
U^+(w)&:=&\{x_{i_1}(t_1)x_{i_2}(t_2)\dotsc x_{i_m}(t_m)\ |\ t_i\in\R_{>0}\},\\
U^-(w)&:=&\{y_{i_1}(t_1)y_{i_2}(t_2)\dotsc y_{i_m}(t_m)\ |\ t_i\in\R_{>0}\}.
\end{eqnarray*} 
These sets are independent of the chosen 
reduced expression and give
\begin{eqnarray*}
U^+(w)&=&U^+_{\ge 0}\cap B^-\dot w B^-,\\
U^-(w)&=& U^-_{\ge 0}\cap B^+\dot w B^+.
\end{eqnarray*}
In particular $U^+_{\ge 0}=\bigsqcup_{w\in W} U^+(w)$ and $U^-_{\ge 0}=
\bigsqcup_{w\in W}U^-(w)$. Moreover $U^+(w)$ and $U^-(w)$ are isomorphic to 
$\R_{>0}^{\ell(w)}$ using the $t_i$ as coordinates. 

Suppose $v\le w$ then $U^+(v)$ can be obtained from $U^+(w)$ by 
letting certain of the $t_i$ coordinates tend to zero. So $U^+(v)$ lies
in the closure of $U^+(w)$. Moreover the condition $v\le w$ is necessary by
the analogous property of the Bruhat decomposition. The same goes for 
the $U^-(w)$. Since $U^+_{\ge 0}$ and 
$U^-_{\ge 0}$ are closed in $G$ by \cite[Proposition~4.2]{Lus:TotPos94}, 
we have
\begin{equation*}
\overline{U^+(w)}=\bigsqcup_{v\le w} U^+(v),\qquad 
\overline{U^-(w)}=\bigsqcup_{v\le w} U^-(w).
\end{equation*} 
Note that in particular 
$\overline{U^+(w_0)}=U^+_{\ge 0}$ and $\overline{U^-(w_0)}=
U^-_{\ge 0}$. The totally positive parts for $U^+$ and $U^-$ 
are defined by 
\begin{equation*}
U^+_{>0}:=U^+(w_0),\qquad U^-_{>0}:=U^-(w_0).
\end{equation*}

\subsection{}\label{s:symmetry}
The totally positive and totally nonnegative parts of the 
flag variety $\mathcal B$ are defined by 
\begin{eqnarray*}
\mathcal B_{>0}&:=&\{y\cdot B^+\ |\ y\in U^-_{>0}\},\\
\mathcal B_{\ge 0}&:=&\overline{\mathcal B_{>0}}.
\end{eqnarray*}
By \cite[Theorem~8.7]{Lus:TotPos94} $\mathcal B_{>0}$ can be
described in a symmetric way as
\begin{equation*}
\mathcal B_{>0}=\{x\cdot B^-\ |\ x\in U^+_{>0} \}.
\end{equation*}
In other words  $\mathcal B_{>0}$ is 
invariant under the automorphism of $G$ (and hence $\mathcal B$) 
which swaps the $x_i(t)$ and the $y_i(t)$. 

\subsection{}\label{s:ansatz} 
The set $\mathcal B_{\ge 0}$ again has a cell decomposition, 
which was conjectured in 
\cite{Lus:TotPos94} and proved in \cite{Rie:CelDec}.
This result was also proved again in 
\cite{MarRie:ansatz} and with explicit descriptions of the cells
given.  We recall the construction from \cite{MarRie:ansatz} below. 

Let $v\le w$ and let $\mathbf w=(i_1,\dotsc, i_m)$ encode a
reduced expression $s_{i_1}\dotsc s_{i_m}$ for $w$. Then there
exists a unique subexpression $s_{i_{j_1}}\dotsc s_{i_{j_k}}$ 
for $v$ in $\mathbf w$ with the 
property that, for $l=1,\dotsc,k$,
\begin{equation*}
s_{i_{j_1}}\dotsc s_{i_{j_{l}}}s_{i_r}>
s_{i_{j_1}}\dotsc s_{i_{j_{l}}}  \quad \text{whenever $j_l< r\le j_{l+1}$,}
\end{equation*}
where $j_{k+1}:=m$. It is the rightmost reduced subexpression for $v$ in 
$\mathbf w$ and we denote it by $\mathbf v=(j_1,\dotsc, j_k)$. 

Then we define
\begin{equation*}
\mathcal R_{v,w}^{>0}:=\left\{g_1\dotsc g_m\cdot B^+\ \left |\
\text{ where $g_r=$}\begin{cases}
\dot s_{i_r}, &\text{ if $r\in\{j_1,\dotsc,j_k\}$,}\\
y_{i_r}(t_r),\ t_r\in\R_{>0} & \text{ otherwise.}
\end{cases}  
\right .\right\}
\end{equation*}
By \cite[Theorem~11.3]{MarRie:ansatz} we have that this definition is independent
of the reduced expression for $w$, and
\begin{equation*}
\mathcal R^{>0}_{v,w}=\mathcal R_{v,w}\cap \mathcal B_{\ge 0}.
\end{equation*}
Moreover the $\mathcal R_{v,w}^{>0}$ are isomorphic to 
$\R_{>0}^{\ell(w)-\ell(v)}$. So this gives an explicit 
decomposition of $\mathcal B_{\ge 0}$ into cells.  

We also write $(B^+\dot w\cdot B^+)_{\ge 0}$ for the intersection 
of the Bruhat cell with $\mathcal B_{\ge 0}$. Then of course
\begin{equation*}
(B^+\dot w\cdot B^+)_{\ge 0}=\bigsqcup_{v;\ v\le w}\mathcal R_{v,w}^{>0}.
\end{equation*}

If $v=1$ then we have $\mathcal R_{1,w}=U^-(w)\cdot B^+$, 
and the above
decomposition of $\mathcal B_{\ge 0}$ 
extends Lusztig's cell decomposition of 
$U^-_{\ge 0}\cong U^-_{\ge 0}\cdot B^+$.

\subsection{}\label{s:symmetry2}
We can apply the symmetry from Section~\ref{s:symmetry} to 
get an alternate description for the $\mathcal R^{>0}_{v,w}$. Namely
let 
\begin{equation*}
\tilde{\mathcal R}_{v,w}^{>0}:=\left\{g_1\dotsc g_m\cdot B^-\ \left |\
\text{ where $g_r=$}\begin{cases}
\dot s_{i_r}\inv, &\text{ if $r\in\{j_1,\dotsc,j_k\}$,}\\
x_{i_r}(t_r),\ t_r\in\R_{>0} & \text{ otherwise.}
\end{cases}  
\right .\right\}
\end{equation*} 
with the same notation as in Section~\ref{s:ansatz}. Then it follows that 
\begin{equation*} 
\tilde{\mathcal R}^{>0}_{v,w}=
\left(B^+\dot v\cdot B^-\cap B^-\dot w\cdot B^-\right)\cap \mathcal B_{\ge 0}\\
=\mathcal R^{>0}_{w w_0, v w_0}.
\end{equation*}

\subsection{}\label{s:reduction}
Suppose $w,w_1,w_2\in W$ with $w=w_1w_2$ and such that the lengths add,
$\ell(w)=\ell(w_1)+\ell(w_2)$. Then there is a well defined map
\begin{eqnarray*}
\pi^w_{w_1}: B^+\dot w\cdot B^+& \to & B^+\dot w_1\cdot B^+\\
B=z\dot w\cdot B^+&\mapsto & z\dot w_1\cdot B^+
\end{eqnarray*} 
where $z\in U^+$. We call $\pi^{w}_{w_1}$ a reduction map.
The element $\pi^{w}_{w_1}(B)$ is uniquely determined by the property
\begin{equation*}
B^+\overset {w_1} \To \pi^{w}_{w_1}(B)\overset{w_2}\To B. 
\end{equation*}
It was proved in \cite{Rie:CelDec} that $\pi^{w}_{w_1}$ preserves 
total positivity. That is, if $B=z\dot w\cdot B^+\in
\mathcal B_{\ge 0}$
then $\pi^{w}_{w_1}(B)$ lies in $\mathcal B_{\ge 0}$. 

Now let $B=g_1\dotsc g_m\cdot B^+\in \mathcal R_{v,w}^{>0}$
where the factors $g_i$ and all the notation are as in 
Section~\ref{s:ansatz}. And additionally suppose 
the reduced
expression for $w$ is a product of reduced expressions for $w_1$ and $w_2$, 
so $w_1=s_{i_1}\dotsc s_{i_{m'}}$ and $w_2=s_{i_{m'+1}}\dotsc s_{i_m}$. 
Then we have explicitly 
\begin{equation*}
\pi^{w}_{w_1}(g_1 g_2\dotsc g_m\cdot B^+)=g_1\dotsc g_{m'}\cdot B^+.
\end{equation*}  
In particular the reduction map restricts to 
a map $\pi^{w}_{w_1}:\mathcal R^{>0}_{v,w}\to \mathcal R^{>0}_{v_{(m')},w_1}$
where $v_{(m')}=s_{i_{j_1}}\dotsc s_{i_{j_{p}}}$ for 
$j_p\le m'< j_{p+1}$.

\section{Closure relations for the cells in $\mathcal B_{\ge 0}$}

The aim of this section is to 
prove the following theorem.

\begin{thm}\label{t:B} Let $v,w\in W$. Then 
$$
\overline{\mathcal R_{v,w}^{>0}}=\bigsqcup_{ 
v\le v'\le w'\le w} 
\mathcal R^{>0}_{v',w'}.
$$ 
\end{thm}

We begin with an easy special case.  
\begin{lem}\label{l:w0} For $v\in W$, 
\begin{equation*}
\overline{\mathcal R^{>0}_{v,w_0}}\cap B^+\dot w_0\cdot B^+=
\bigsqcup_{\ v'\ge v}\mathcal R^{>0}_{v',w_0}.
\end{equation*}
\end{lem}
\begin{proof}
By Section~\ref{s:symmetry2} we can rewrite 
\begin{equation*}
\mathcal R_{v',w_0}^{>0}=\tilde{\mathcal R}^{>0}_{1,v' w_0}=
U^+(v' w_0)\dot w_0\cdot B^+.
\end{equation*}
Now $v'\ge v$ implies $v'w_0\inv\le vw_0$. Hence by
Section~\ref{s:totpos} we have the inclusion $U^+(v'w_0)\subseteq 
\overline{U^+(vw_0 )}$, and therefore 
$\mathcal R_{v',w_0}^{>0}\subseteq\overline {\mathcal R_{v,w_0}^{>0}}.$ 
The opposite inclusion follows from the closure relations of 
Bruhat decomposition.
\end{proof}

\begin{lem}\label{l:phi}
There is a homeomorphism
\begin{equation*}
\phi=\phi_w:U^-(w_0w\inv)\x (B^+ \dot w\cdot B^+)_{\ge 0}
\overset\sim\To
\bigsqcup_{v\le w}\mathcal R_{v, w_0}^{>0}
\end{equation*}
defined by $\phi(u, B):=u\cdot B$. 
Moreover 
\begin{equation*} 
\phi_w(U^-(w_0 w\inv)\x \mathcal R^{>0}_{v,w})=
\mathcal R^{>0}_{v,w_0}.
\end{equation*}
\end{lem}

\begin{proof}
Choose a reduced expression $\mathbf w_0=(i_1,\dotsc, i_N)$
for $w_0$ such that $(i_1,\dotsc i_r)$ is a reduced expression
for $w_0w\inv$.
Then using the parameterizations of
$U^-(w_0 w\inv)$ and $\mathcal R^{>0}_{v,w_0}$ described in  
Sections~\ref{s:totpos} and \ref{s:ansatz}, respectively,
we see that 
\begin{eqnarray*}
U^-(w_0 w\inv)\x\mathcal R^{>0}_{v,w}&\overset\sim\to &
\mathcal R^{>0}_{v,w_0},\\
(u,B)&\mapsto & u\cdot B.
\end{eqnarray*}
is an isomorphism. Explicitly we have 
\begin{equation*}
\mathcal R_{v,w_0}=\left\{
y_{i_1}(t_1)\dotsc y_{i_r}(t_r)g_{r+1}\dotsc g_{N}\cdot B^+\ \left |
\ g_{r+1}\dotsc g_{N}\cdot B^+\in\mathcal R_{v,w}^{>0}
\right .\right\}
\end{equation*} 
Now applying the reduction map $\pi^{w_0}_{w_0w\inv}$
from Section~\ref{s:reduction} we get
\begin{eqnarray*}
\mathcal R_{v,w_0}^{>0}&\To &
\mathcal R^{>0}_{1,w_0w\inv}\\
y_{i_1}(t_1)\dotsc y_{i_r}(t_r)g_{r+1}\dotsc g_{N}\cdot B^+
&\mapsto &y_{i_1}(t_1)\dotsc y_{i_r}(t_r)\cdot B^+
\end{eqnarray*}
Note that $\pi^{w_0}_{w_0w\inv}$ is defined on the whole
Bruhat cell  $B^+\dot w_0\cdot B^+$. We can therefore 
combine these maps for varying $v$ and compose with the isomorphism  
$U^-\cdot B^+ \overset\sim\To U^-$, to get 
\begin{equation*}
p_1: \bigsqcup_{v\le w} \mathcal R^{>0}_{v,w_0}\to U^-(w_0 w\inv).
\end{equation*}
The inverse to $\phi$ is now given by
\begin{eqnarray*}
\psi\ :\ \  \bigsqcup_{v\le w} \mathcal R^{>0}_{v,w_0}&\to& U^-(w_0 w\inv)\x
(B^+\dot w\cdot B^+)_{\ge 0}\\
  B &\mapsto & \left(p_1(B),p_1(B)\inv\cdot B\right).
\end{eqnarray*}
\end{proof}

\begin{lem}\label{l:w} Let $v\le v'\le w\in W$. Then 
$$
\mathcal R^{>0}_{v',w}\subseteq \overline{\mathcal R_{v,w}^{>0}}.
$$
\end{lem}

\begin{proof}
Using Lemma~\ref{l:w0} we see that 
\begin{equation*}
\mathcal R_{v',w_0}^{>0}
\subset \overline{\mathcal R_{v,w_0}^{>0}}\cap \bigsqcup_{x\le w}
\mathcal R_{x,w_0}^{>0}. 
\end{equation*}
Applying $\psi:=\phi\inv$ from Lemma~\ref{l:phi}
to both sides of this inclusion gives
\begin{equation*}
\psi(\mathcal R_{v',w_0}^{>0})
\subset \psi(\overline{\mathcal R_{v,w_0}^{>0}}\cap \bigsqcup_{x\le w}
\mathcal R_{x,w_0}^{>0})\subset \overline{\psi(\mathcal R^{>0}_{v,w_0})}. 
\end{equation*}
Therefore
\begin{equation*}
U^-(w_0w\inv)\x\mathcal R^{>0}_{v',w}=
\psi(\mathcal R^{>0}_{v',w_0})\subset 
\overline{\psi(\mathcal R^{>0}_{v, w_0})}=\overline{
U^-(w_0 w\inv)\x \mathcal R^{>0}_{v,w}},
\end{equation*}
where we may take the closure on the right hand side 
to be the closure inside the domain of $\phi$, that is
inside $U^-(w_0 w\inv)\x (B^+\dot w\cdot B^+)_{\ge 0}$.  
It follows that $\mathcal R^{>0}_{v',w}\subseteq\overline
{\mathcal R^{>0}_{v,w}}$.
\end{proof}

\begin{proof}[Proof of Theorem~\ref{t:B}]
By \cite{Lus:TotPos94} $\mathcal B_{\ge 0}$ is 
symmetric with respect to interchanging $B^+$ and $B^-$, see 
Section~\ref{s:symmetry}. Therefore it follows that Lemma~\ref{l:w} 
also holds for the $\tilde{\mathcal R}^{>0}_{v,w}$ defined in 
Section~\ref{s:symmetry2}. That is,
$
\tilde {\mathcal R}^{>0}_{v',w}\subseteq 
\overline{\tilde{\mathcal  R}^{>0}_{v,w}}
$  
whenever $v\le v'\le w$. 
Now if $v\le v'\le w'\le w$, then 
$$
\mathcal R^{>0}_{v',w'}\subseteq
\overline{\mathcal R^{>0}_{v',w}}=
\overline{\tilde{\mathcal R}^{>0}_{ww_0,v'w_0}}
\subseteq\overline{\tilde{\mathcal R}^{>0}_{ww_0,vw_0}}=
\overline{\mathcal R^{>0}_{v,w}}.
$$
This shows the inclusion $\supseteq$ in the statement of Theorem~\ref{t:B}. 

The other inclusion is clear from the closure relations of 
Bruhat decomposition. $\mathcal R^{>0}_{v',w'}\cap
\overline{\mathcal R^{>0}_{v,w}}\ne \emptyset$ implies on the one hand
$\mathcal R_{v',w'}^{>0}\cap  \overline{B^-\dot v\cdot B^+}\ne\emptyset$ and 
on the other hand $\mathcal R_{v',w'}^{>0}\cap
\overline{ B^+\dot w \cdot B^+}\ne \emptyset$.  
So we must have $v\le v'$ and $w'\le w$. 
\end{proof}

We note that Theorem~\ref{t:B} 
implies  
$\overline{\mathcal R_{v,w}}\cap 
\mathcal B_{\ge 0}=\overline{\mathcal R^{>0}_{v,w}}$.

\section{Lusztig's decomposition of $\mathcal P^J$}\label{s:strata}

The stratification of $\mathcal B$ into smooth pieces $\mathcal R_{v,w}$ has
an analogue for partial flag varieties introduced by Lusztig in 
\cite{Lus:PartFlag}. 

Consider a triple of Weyl group elements $x,u,w\in W$ with $x\in W^J_{max}$,
$w\in W^J$ and $u\in W_J$. Then 
$\mathcal P^J_{x,u,w}\subset \mathcal P^J$ is defined as the 
set of all $P\in \mathcal P^J$ such that there exist Borel subgroups 
$B_L$ and $B_R$ inside $P$ satisfying
\begin{equation*}
B^+\overset w\To B_L\overset u\To B_R\overset {x\inv w_0}\To B^-.
\end{equation*}
An equivalent characterization of $\mathcal P^J_{x,u,w}$ is 
\begin{equation*}
\mathcal P^J_{x,u,w}=\pi^J(\mathcal R_{x,w u})=\pi^J(\mathcal R_{x u\inv, w}).
\end{equation*}
It is not hard to see that $B_L$ and $B_R$ are uniquely determined
as the Borel subgroups in $P$ `closest to' $B^+$ respectively $B^-$
with regard to their relative position, and the projection maps 
$\mathcal R_{x,wu}\to \mathcal P^J_{x,u,v}$ and $\mathcal R_{x u\inv,w}\to
\mathcal P^J_{x,u,v}$ are isomorphisms. 
In particular $\mathcal P^J_{x,u,w}$ is nonempty if and only
if $x\le w u$, in which case it is smooth of dimension 
$\ell(w)+\ell(u)-\ell(x)$.  

Let us denote the indexing set for this decomposition of $\mathcal P^J$
by $Q^J$. So
\begin{equation*}
Q^J:=\{(x,u,w)\in W^J_{max}\x W_J\x W^J\ |\ x\le wu \}.
\end{equation*} 

We define the following partial oder on $Q^J$.
\begin{defn}\label{d:poset}
Let $(x',u',w')$ and $(x,u,w)$ in $Q^J$. Then define
\begin{equation*}
(x',u',w')\le (x,u,w)
\end{equation*}
if and only if there exist $u'_1,u'_2\in  W_J$ satisfying 
$u'_1 u'_2=u'$ with $\ell(u'_1)+\ell(u'_2)=\ell(u')$, and such that 
\begin{equation}\label{e:inequalities}
x u\inv\le x'{u_2'}\inv\le w'u'_1\le w.
\end{equation} 
\end{defn}

\section{Totally nonnegative cells in 
$\mathcal P^J$ and their closure relations}\label{s:P}

The totally positive and nonnegative parts of 
$\mathcal P^J$ are defined in
\cite{Lus:IntroTotPos} by
\begin{eqnarray*}
\mathcal P^J_{>0}&=&\pi^J(\mathcal B_{>0}),\\
\mathcal P^J_{\ge 0}&=&\pi^J(\mathcal B_{\ge 0}).
\end{eqnarray*} 
Since $\pi^J$ is closed it follows that $\mathcal P^J_{\ge 0}=
\overline{\mathcal P^J_{>0}}$.

We decompose $\mathcal P^J_{\ge 0}$ by intersecting it with the 
strata $\mathcal P^J_{x,u,w}$ from Section~\ref{s:strata}. 
From the definitions and the fact that reduction preserves total
positivity it follows that (\cite[Lemma 3.2]{Rie:PhD}) 
\begin{equation*}
\mathcal P^J_{x,u,w;>0}:=\mathcal P^J_{x,u,w}\cap\mathcal P^J_{\ge 0}
=\pi^J(\mathcal R_{x,wu}^{>0})=\pi^J(\mathcal R_{x u\inv, w}^{>0}).
\end{equation*}
Keeping in mind that 
$\pi^J:\mathcal R_{x,wu}\to\mathcal P^J_{x,u,w}$, say, 
is an isomorphism, we see that 
\begin{equation*}
\mathcal P^{J}_{x,u,w;>0}\cong \mathcal R^{>0}_{x,wu}\cong
\mathcal \R_{>0}^{\ell(w)+\ell(u)-\ell(x)},
\end{equation*}
for any triple $(x,u,w)\in Q^J$.

We will prove the following theorem describing 
the closure relations between the strata $\mathcal P^{J}_{x,u,w;>0}$
in $\mathcal P^J_{\ge 0}$.
 
\begin{thm}\label{t:P} Let $(x,u,w)\in Q^J$ with partial
order as in Definition~\ref{d:poset}, then
\begin{equation*}
\overline{\mathcal P^J_{x,u,w;>0}}=
\bigsqcup_{(x',u',w')\le (x,u,w)}\mathcal P^J_{x',u',w';>0}.
\end{equation*}
\end{thm}

\begin{proof} Note first that we have
$\overline{\mathcal P^J_{x,u,w;>0}}=
\pi^J(\overline{\mathcal R^{>0}_{x u\inv, w}})=
\pi^J(\overline{\mathcal R^{>0}_{x,w u}})$.

Now suppose $(x',u',w')\le (x,u,w)$. So we have 
$u'_1$ and $u'_2$ as in Definition~\ref{d:poset}. Let 
$P'\in\mathcal P^{>0}_{x',u',w'}$ with its associated 
Borel subgroups $B'_L\in \mathcal R_{x'{u'}\inv,w'}^{>0}$ 
and $B'_R\in\mathcal R^{>0}_{x',w'u'}$ such that $B'_L\subset P'$
and $B'_R\subset P'$. In particular  
\begin{equation*}
B^+\overset{w'}\To B'_L\overset{u'}\To B'_R\overset {{x'}\inv w_0}\To B^-.
\end{equation*}
Let $B':=\pi^{w'u'}_{w'u'_1}(B'_R)$ using the reduction map from 
Section~\ref{s:reduction}. Then we have
\begin{equation*}
B^+\overset{w'}\To B'_L\overset{u'_1}\To B'\overset{u'_2}\To B'_R
\overset{{x'}\inv w_0}\To B^-.
\end{equation*}
Since reduction preserves total positivity we have 
$B'\in \mathcal R^{>0}_{x'{u'_2}^{-1},w' u'_1}$. Also $\pi^J(B'_L)= P'$ and
the fact that $(B'_L,B')$ have relative position $u'_1\in W_J$ 
implies that $\pi(B')=P'$. Now by \eqref{e:inequalities} together
with Theorem~\ref{t:B} it follows that   
$\mathcal R^{>0}_{x'{u'_2}^{-1},w'u'_1}\subset
\overline{\mathcal R^{>0}_{x u^{-1},w}}$. Therefore
$P'=\pi (B')\in \pi(\overline{\mathcal R^{>0}_{x u^{-1},w}})=\overline 
{\mathcal P^J_{x,u,w;>0}}$. This proves 
the inclusion $\supseteq$.

For the opposite inclusion suppose 
that $\mathcal P_{x',u',w';>0}\cap\overline 
{\mathcal P_{x,u,w;>0}}\ne\emptyset$. So let $P'$ be an element
of this intersection.  
Then there exists a  
$\tilde B'\in\overline{\mathcal R^{>0}_{xu^{-1},w}}$ such that 
$\pi(\tilde B')=P'$. Since $\tilde B'\subset P'$ we have  
$u_1,u_2\in W_J$ such that 
$$
\begin{matrix}
B^+&\overset{w'}\longrightarrow &B'_L & &\overset{u'}
\longrightarrow & &B'_R&\overset{(x')^{-1}w_0}\longrightarrow&B^-\\
&&\ u_1&\searrow& &\nearrow &u_2\ &&\\
&&&   & \tilde B' &  &&&
\end{matrix}.
$$
If it happens to be the case that $\ell(u_1 u_2)=\ell(u_1)+\ell(u_2)$ then
$u'=u_1 u_2$. So we can set $u'_1=u_1$ and $u'_2=u_2$ and are done. 
Otherwise there exists a simple
reflection $s=s_{i_1}\in W_J$ such that $u_1 s\le u_1$ and 
$s u_2\le u_2$. We have 
$$
\begin{matrix}
B'_L   &        && & \overset{u'} \longrightarrow &&& &B'_R\\
\ \   u_1 s&\searrow& &     &&                        &&\nearrow &s u_2\ \ \\
        &      &B'_{L,1}&   &  
\overset{\text{$s$ or $1$}}\longrightarrow     &     &B'_{R,1} &  &\\
        &      &\ \   s &\searrow&   &\nearrow & s\ \ & & \\
& &   && \tilde B'& &  & &.
\end{matrix}
$$
Here $B'_{L,1}$ is obtained from $\tilde B'$ by reduction,
$B'_{L,1}=\pi^{wu_1}_{wu_1s}(\tilde B')$. This implies
$B'_{L,1}\in\mathcal B_{\ge 0}$ (see Section~\ref{s:reduction}).
Using the inequalities
$$ x u^{-1}\le x' u_2^{-1}\le x'{(su_2)^{-1}}, 
\qquad w'u_1 s\le w'u_1\le w,  
$$
and Theorem~\ref{t:B} it follows that
$$
B'_{L,1}\in\mathcal R^{>0}_{x' u_2^{-1},w' u_1 s}\sqcup 
\mathcal R^{>0}_{x'(su_2)^{-1},w'u_1 s}\subseteq 
\overline{\mathcal R^{>0}_{xu^{-1},w}}.
$$
Therefore we are now in the analogous situation as before
$$
\begin{matrix}
B^+&\overset{w'}\longrightarrow &B'_L & &\overset{u'}
\longrightarrow & &B'_R&\overset{(x')^{-1}w_0}\longrightarrow&B^-\\
&&\ u_1^{(1)}&\searrow& &\nearrow &u_2^{(1)}\ &&\\
&&&   & B'_{L,1} &  &&&
\end{matrix}.
$$
with totally nonnegative $B'_{L,1}$, 
but where $u_1^{(1)}=u_1 s< u_1$. If again
$\ell(u_1^{(1)}u_2^{(1)})\ne \ell(u_1^{(1)})+\ell(u_2^{(1)})$ 
then we can 
repeat the argument above. So we replace $u_1^{(1)}$
with a shorter
$u_1^{(2)}=u_1^{(1)}s_{i_2}$, and $B'_{L,1}$ with $B'_{L,2}$ also
in $\overline{\mathcal R^{>0}_{xu^{-1},w}}$. Iterating this 
process we must eventually arrive at a case where 
$\ell(u_1^{(k)}u_2^{(k)})=\ell(u_1^{(k)})+\ell(u_2^{(k)})$ 
(after at most $k=\ell(u_1)$ steps), at which point we set
$u_1^{(k)}=u_1'$ and $u_2^{(k)}=u_2'$, and
$B'_{L,k}\in\overline{\mathcal R^{>0}_{xu^{-1},w}}$ implies
the inequalities \eqref{e:inequalities}. 
\end{proof}

\section{Closure relations for the strata $\mathcal P^J_{x,u,w}$ 
in $\mathcal P^J$}

Let $\mathbb K$ be again as in Section~\ref{s:prelim}. 
Then we have a decomposition  
\begin{equation*}
\mathcal P^J=\bigsqcup_{(x,u,w)\in Q^J}
\mathcal P^J_{x,u,w}
\end{equation*}
defined as in Section~\ref{s:P}.  
We can now deduce the analogue of
Theorem~\ref{t:P} for the strata $\mathcal P^J_{x,u,w}$
in $\mathcal P^J$.

Again the full flag variety case needs to be treated first. 
After that the proof proceeds in the same way as for Theorem~\ref{t:P}.

\begin{prop}\label{p:B} Let $v,w\in W $ with  $v\le w$. Then
\begin{equation}
\overline{\mathcal R_{v,w}}=\bigsqcup_{v\le v'\le w'\le w}\mathcal R_{v',w'}. 
\end{equation}
\end{prop}
The above result does not seem to appear in the literature, so we include 
a quick proof. It is however well-known to experts, 
\cite{Lus:PersonalCommunication05}.   
\begin{proof}
The inclusion $\subseteq$ follows from the closure relations for Bruhat
decomposition.

Let us prove the inclusion $\supseteq$. 
It suffices by $B^+$-$B^-$ symmetry to show that 
$\mathcal R_{v',w}\subseteq \overline{\mathcal R_{v,w}}$ whenever
$v\le v'$ (as in the proof of Theorem~\ref{t:B}). 
Consider the map 
\begin{eqnarray*}
\gamma~: B^+\dot w\cdot B^+ \x (U^-\cap \dot w U^-\dot w\inv )
&\To & \dot w U^-\cdot B^+,\\
(B,y) &\mapsto & y\cdot B.
\end{eqnarray*}
It is an isomorphism and by restriction 
gives rise to isomorphisms
\begin{equation*}
\gamma_v~: 
\mathcal R_{v,w}\x (U^-\cap \dot w U^-\dot w\inv)\To 
B^-\dot v\cdot B^+\ \cap\ \dot w U^-\cdot B^+, 
\end{equation*}
as in \cite[Section 1.4]{Lus:IntroTotPos}. Now if $v'\ge v$ then 
we have 
\begin{equation*}
B^-\dot v'\cdot B^+\cap \dot w U^-\cdot B^+\subset
\overline{B^-\dot v\cdot B^+}\cap \dot w U^-\cdot B^+.
\end{equation*}
Applying $\gamma\inv$ to this inclusion we see that
$\mathcal R_{v',w}\subset \overline{\mathcal R_{v,w}}$ and
the proposition follows.
\end{proof}

Now consider the partial flag variety case. Since $\pi^J$ is proper
we have that $\overline{\mathcal P^J_{x,u,w}}=
\pi^J(\overline{\mathcal R_{x,wu}})=
\pi^J(\overline{\mathcal R_{xu\inv, w}})$. 
The following proposition follows from the same proof 
as Theorem~\ref{t:P}, only using Proposition~\ref{p:B} 
in place of Theorem~\ref{t:B} and leaving 
out the positivity considerations.

\begin{prop}\label{p:P} Let $(x,u,w)\in Q^J$ with partial
order as in Definition~\ref{d:poset}, then
\begin{equation*}
\overline{\mathcal P^J_{x,u,w}}=
\bigsqcup_{(x',u',w')\le (x,u,w)}\mathcal P^J_{x',u',w'}.
\end{equation*}
\qed
\end{prop} 

Note that the varieties $\mathcal P_{x,u,w}$ can be defined over any field $\mathbb K$ and are $\mathbb K$-rational. Assuming that $\mathbb K$ is infinite, it follows that Propositions~\ref{p:B} and \ref{p:P} hold also in this setting.


\def\cprime{$'$}
\providecommand{\bysame}{\leavevmode\hbox to3em{\hrulefill}\thinspace}
\providecommand{\href}[2]{#2}

\end{document}